\documentclass[%
 aip,
 amsmath,amssymb,
 reprint,%
]{revtex4-2}
 \usepackage{hyperref} 
\usepackage{graphicx}
\usepackage{caption}
\usepackage{textgreek}
\usepackage{dcolumn}

\usepackage{nccmath}

\usepackage{amsmath}
\usepackage{subcaption}
\DeclareFontFamily{OMX}{cmexb}{}
\DeclareFontShape{OMX}{cmexb}{b}{n}{<->cmexb10}{}
\DeclareSymbolFont{boldoperators}{OT1}{\familydefault}{b}{n}
\DeclareSymbolFont{boldlargesymbols}{OMX}{cmexb}{b}{n}
\DeclareMathDelimiter{\lbrackbold}{\mathopen} {boldoperators}{"5B}{boldlargesymbols}{"02}
\DeclareMathDelimiter{\rbrackbold}{\mathclose}{boldoperators}{"5D}{boldlargesymbols}{"03}
\usepackage{verbatim}
\makeatletter
\newenvironment{boldbmatrix}{\left\lbrackbold\env@matrix}{\endmatrix\right\rbrackbold}
\makeatother
\usepackage{bm}
\usepackage{mathtools}
\usepackage{mleftright}
\usepackage{xparse}
\usepackage{textcomp} 
\newcommand{\qdist}[1]{\ifmmode\langle#1\rangle\else\textlangle#1\textrangle\fi}

\begin{document}

\title[]{On the Defense of the Recent Solution to Hilbert's Sixteenth Problem: Clarifying Misinterpretations and the Incorrect Conclusions in Buzzi and Novaes's note}




\author{$^1$Vinícius Barros da Silva} 
\email{vinicius.barros@unesp.br}
\author{$^2$João Peres Vieira and $^1$Edson Denis Leonel}
\address{$^1$Department of Physics, Universidade Estadual Paulista ``J{\'u}lio de Mesquita 
Filho", Campus de Rio Claro, S{\~a}o Paulo 13506-900, Brazil\\
$^2$Department of Mathematics, Universidade Estadual Paulista ``J{\'u}lio de Mesquita 
Filho", Campus de Rio Claro, S{\~a}o Paulo 13506-900, Brazil
} 


\date{\today}

\begin{abstract}
Recently, the covariant formulation of the geometric bifurcation theory, developed in a previous paper, has been applied to two elementary problems- the study of limit cycles of dynamical systems and the second part of Hilbert's sixteenth problem. First, it has been shown that dynamical systems with more than one limit cycle are understood to be those in which the scalar curvature R is positive and its magnitude diverges to infinity at different singular points. In the second, it has been demonstrated that $n$th-degree polynomial systems have the maximum number of $2(n-1)\bm{(}4(n-1)-2\bm{)}$ limit cycles with their relative positions determined by the singularities of the magnitude of R, thus providing a successful response to the original Hilbert's challenge. It is the purpose of this letter to point out that Buzzi and Novaes's note is incorrect and leads to erroneous results as an immediate consequence of their superficial reading of our work, the omission of critical aspects of the GBT to induce doubt, and reliance on incorrect assumptions or interpretations that deviate from those established in the framework of GBT.

\end{abstract}

\maketitle


%

\section{\label{sec:Intro}Introduction}

The covariant formulation of the geometric bifurcation theory (GBT) has been applied in two important applications in the third article of this series\cite{Vinicius,Viniciusprime,Viniciussecond}. These applications were the study of isolated periodic trajectories of dynamical systems governed by differential equations and the second part of Hilbert's sixteenth problem, questioning the maximum number and position of Poincaré's limit cycles for polynomial differential equations of degree $n$.

As a result of several investigations of systems with different numbers of limit cycles, we found that dynamical systems with more than one limit cycle are understood to be those in which the scalar curvature R is positive and its magnitude diverges to infinity at different singular points. Consequently, the total number of distinct divergences of $|\mbox{R}|$ to infinity provides the maximum number of limit cycles. In addition, we have demonstrated that $n$th-degree polynomial systems have the maximum number of $2(n-1)\bm{(}4(n-1)-2\bm{)}$ limit cycles with their relative positions determined by the singularities of the magnitude of R.

However, there have been two concerns regarding those results. First, the apparent contradiction between the Hilbert number demonstrated in our original paper and those established in previous studies. Second, the validity of the covariant definition of limit cycles. This, together with groundless theoretical objections, has raised doubts concerning the validity of our solution for Hilbert's sixteenth problem.

In this paper, we will lay those doubts to rest. We begin by demonstrating how the underlying assumptions that lead to the logarithmic growth of the Hilbert number deviate from the original problem and other established results in the literature. Finally,  we conclude by showing the validity of the covariant definition of limit cycles and by pointing out that the erroneous conclusions of Buzzi and Novaes arise from a superficial reading of Ref.~\cite{Vinicius}, the omission of critical aspects of the GBT in their analysis, and reliance on incorrect assumptions or interpretations that deviate from those of our framework.


\section{\label{sec:GT} The Hilbert number}

In order to understand the apparent contradiction between the Hilbert number, $\mbox{H}_{n}$, demonstrated in our original paper\cite{Vinicius} and those suggested throughout history, it is important to recognize first the nature of Hilbert’s sixteenth problem, the usual methods applied to demonstrate that $\mbox{H}_{n}$ increases with $n^{2} \,\text{log}\,n$  and the underlying hypotheses, mostly omitted, that support this latter result.

In general, mathematical challenges can be naturally divided into two distinct categories: determination and demonstration problems. According to Pólya and Conway, a determination problem aims to find a certain object—the unknown of the problem—using the unknown, the data, and the condition. Conversely, the aim of a demonstration problem is to show conclusively that a certain clearly stated assertion is true, or else to show that it is false by working with a hypothesis and a conclusion~\cite{Polya}.

In other words, the division into determination and demonstration problems highlights the distinct ways mathematical reasoning operates. Determination problems often involve creative exploration to derive a solution from given conditions, while demonstration problems focus on rigorous logical proof to validate or refute specific claims. The structure of a problem significantly shapes the methods used to address it and, of course, the conclusions drawn.

Bearing this in mind, let us now turn our attention to the statements of the second part of Hilbert’s Sixteenth Problem~\cite{Hilbert}: 
\begin{quote}
~~~``In connection with this purely algebraic problem, I wish to bring forward a question which, it seems to me, may be attacked by the
same method of continuous variations of coefficients, and whose answer is of corresponding value for the topology of families of curves defined by differential equations. This is the question as to maximum number and position of Poincaré's boundary cycles (limit cycles) for a differential equation of the first order and degree of the form
\begin{equation*}
\frac{dy}{dx}=\frac{\mbox{Y}}{\mbox{X}},
\end{equation*}
where $\mbox{X}$ and $\mbox{Y}$ are rational integral functions of the $n$th degree in $x$ and $y$ ".
\end{quote}

As one may observe, Hilbert's problem is fundamentally a determination problem, as it seeks to identify the maximum number of limit cycles for a polynomial system of differential equations of degree no greater than $n$. That is to say, Hilbert does not introduce any conditional relationships or auxiliary assumptions. He explicitly frames the problem within the context of the given data, allowing no room for additional hypotheses about the nature of the polynomials or their degree\cite{Gaiko}. This clarity leaves no basis for imposing further constraints or supplementary conditions.

Having outlined the nature of Hilbert's challenge, let us now turn to the predictions upon the supposed logarithmic growth of $\mbox{H}_{n}$.


The historical results that suggest that the Hilbert number must grow as $n^{2} \,\text{log}\,n$ relie heavily on perturbation methods and assume that the degree of the polynomials can be doubled with each perturbation\cite{Chris,Li}. According to Christopher and Lloyd~\cite{Chris}:
\begin{quote}
~~~``The system is transformed by a singular transformation which quadruples the number of limit cycles, while doubling the degree of the system.”
\end{quote}

As we may observe, their results are grounded in a particular hypothesis concerning polynomial degrees, specifically the modification of polynomial degrees—a premise that Hilbert did not establish in his original formulation of the problem. While Christopher and Lloyd's approach is valuable for constructing systems with a high number of limit cycles, it fails to accurately reflect the behavior of the unperturbed system. Moreover, Li, Chan, and Chung observed in 2002 that~\cite{Li}:
\begin{quote}
``Recently, Christopher and Lloyd~[6] showed that $H(2^{k}-1) \geq 4^{k-1}(2k-\frac{35}{6})+3.2^{k}-\frac{5}{3}$(for example $H(7) \geq 25$) by perturbing some families of closed orbits of a Hamiltonian system sequence in small neighbourhoods of some center points and using a ``quadruple transformation". The method given by [6] is very interesting. Unfortunately, the computation of a lower bound is not correct".
\end{quote}


This statement underscores a critical flaw in Christopher and Lloyd’s approach, highlighting that the prediction upon which $\mbox{H}_{n}$ grows as $n^{2} \,\text{log}\,n$ is incorrect. 

Furthermore, we must remark on two important aspects of the approach of  Li, Chan, and Chung, which most researchers forget to remember or remember to forget. First, while Li, Chan, and Chung have corrected the numerical results obtained by Christopher and Lloyd, their work still operates under the same premise—that the degree of the polynomials can be modified using perturbation methods. That continued reliance on an unsupported hypothesis remains inconsistent with the original determination problem outlined by Hilbert. Second, their method is devoted to the study of the weakened version of Hilbert’s 16th problem as proposed by V.I. Arnold in 1977. This version focuses exclusively on the analysis of symmetric planar polynomial Hamiltonian systems, rather than general polynomial differential equations of degree $n$ as originally queried by Hilbert in 1900.


Surprisingly, previous studies have also demonstrated that $\mbox{H}_{n}$ grows proportionally to $n^{2}$ when perturbation methods are employed to study Hilbert's sixteenth problem without relying on the assumption that the degree of polynomials can be modified, in remarkable agreement with the results obtained by Il'jašenko~\cite{Yas}. Further, this conclusion is consistent with the findings of Sheng, who also demonstrated that Hilbert number indeed scales with $n^{2}$ through the use of algebraic geometry, Bézout's theorem, and polynomial decomposition techniques~\cite{Pingxing}. More precisely, Sheng's analysis established that the maximum number of limit cycles is bounded by $d_{f} \cdot d_{g}$, where $d_{f}$ and $d_{g}$ are the degrees of the polynomials defining the system, thereby reinforcing the quadratic growth when $d_{f}=d_{g}=n$. 


In conclusion, while the work of Christopher and Lloyd provides valuable insights into systems with perturbations, it does not solve Hilbert's original problem. By introducing perturbations, their method modifies the system in ways that obscure the unperturbed behavior. Moreover, although Li, Chan, and Chung addressed errors in Christopher and Lloyd’s results, their approach still relies on the same unsupported hypothesis of modifying the polynomial’s degree. This assumption is inconsistent with Hilbert's original challenge, as the degree of the polynomials cannot exceed  $n$.


On the other hand, our results, derived through the application of the GBT, conclusively demonstrate that the Hilbert number grows proportionally to $n^{2}$.  This finding aligns remarkably well with the results of Il'jašenko by using perturbation methods without assuming the hypothesis that the degree of the polynomials can be modified. In addition, it agrees with the conclusions of the rigorous algebraic approach of Sheng, who corroborated the quadratic growth of the Hilbert number using geometrical methods. Therefore, by maintaining the integrity of Hilbert's original formulation, our work provides a faithful and robust solution to the problem as it was originally posed.

\section{\label{sec:GT} The Covariant definition of limit cycles} 

We now address Buzzi and Novaes' concerns regarding the covariant definition of limit cycles introduced in our recent work~\cite{Vinicius}. They claimed that the positivity of R is not a necessary condition for the existence of periodic trajectories. Buzzi and Novaes reach their erroneous conclusion through a superficial examination of scalar curvature at equilibrium points. Their analysis, along with a few purported ``counterexamples", raised doubts regarding the validity of the covariant definition of limit cycles.

In this section, we shall lay those doubts to rest. It is essential to note that Buzzi and Novaes' observations demonstrate a fundamental misinterpretation of the geometrical methods outlined in Ref.~\cite{Vinicius}. Their conclusions stem from a superficial reading of our work, overlooking critical aspects of the analysis of R. Contrary to their claims, the study and interpretation of R within Riemannian geometry extend beyond mere evaluations at equilibrium points. In fact, the study of R sign is fundamentally linked to the Poincaré index theorem and its geometric counterpart, the Gauss-Bonnet theorem, as rigorously detailed and discussed in Ref.~\cite{Vinicius}. This oversight by Buzzi and Novaes has led to their erroneous conclusions.

In light of the preceding discussion and the appropriate GBT methods, we will now re-examine their alleged ``counterexamples". The analysis to be carried out aligns with the one established in section 03 of Ref.~\cite{Vinicius}, which we will summarize for brevity. Through this examination, we will demonstrate that these examples effectively illustrate the robustness and validity of the GBT approach, reinforcing the accuracy of our findings.\\
 
\subsection{Example 01}
The first supposed ``counterexample" is the following, 
\begin{equation}\label{eq:ex01}
\left\{\begin{array}{lcl} &\displaystyle\frac{ds_{1}}{d\tau}&=\Phi_{1}\left(s_{1},s_{2}\right)=-s_{2} + s_{1}(s_{1}^2 + s_{2}^2 - 1) ,\\ \\ &\displaystyle\frac{ds_{2}}{d\tau}&=\Phi_{2}\left(s_{1},s_{2}\right)= s_{1} + s_{2}(s_{1}^2 + s_{2}^2 - 1) , \end{array}  \right. 
\end{equation}
which $p_{1}=(0,0)$ is an equilibrium point. According to Buzzi and Novaes's analysis, the scalar curvature of Eq.~(\ref{eq:ex01}) is negative near $p_{1}=(0,0)$. However, the sign of R in the GBT is not determined by merely evaluating equilibrium points in R but rather through the Poincaré and the Gauss-Bonnet theorems.

Bearing this in mind, let us now examine the sign of R in the framework of the covariant formulation of GBT. From 	Eq.~(\ref{eq:ex01}), it is not difficult to recognize that it has two varying parameters: $s_{1}$ and $s_{2}$. So the Riemannian metrical structure of the parameter space $X=(s_{1},s_{2})$ corresponds to a two-dimensional Riemannian geometry. Considering this, we work forward from GBT by deriving the invariant positive definite metric. From Eqs.~(5) and (7) of Ref.~\cite{Vinicius} and~(\ref{eq:ex01}), we have
\begin{widetext} 
\begin{equation}\label{metric_eq:ex01} 
 d\ell^{2}=2\pmb{\big[} 
2 \left((3 s_{1}^2+s_{2}^2-1\right)^2+(2 s_{1} s_{2}+1)^2)\,\pmb{\big]}^{2}(ds_{1})^{2} 
+2\pmb{\big[} 2 \left((s_{1}^2+3 s_{2}^2-1)^2+(2 s_{1} s_{2}-1)^2\right)\,\pmb{\big]}(ds_{2})^{2}. 
\end{equation}

The metric of Eq.~(\ref{metric_eq:ex01}) imposes a curvature R on the two-dimensional parameter space $X$ of the system above. From Eqs.~(8) of Ref.~\cite{Vinicius} and (\ref{eq:ex01}), we then encounter 
\begin{multline}\label{eq:R_ex01} 
\mbox{R}= A^{-2}(s_{1},s_{2})\pmb{\Bigg[}+4 (18 s_{1}^{10}-3 s_{1}^8 (18 s_{2}^2+17)-80 s_{1}^7 s_{2}+s_{1}^6 (-764 s_{2}^4+116 s_{2}^2+92)+48 s_{1}^5 (s_{2}^3+s_{2})\\ -2 s_{1}^4 (382 s_{2}^6-295 s_{2}^4+38 s_{2}^2+30)-48 s_{1}^3 s_{2}^5-2 s_{1}^2 (27 s_{2}^8-58 s_{2}^6+38 s_{2}^4+12 s_{2}^2-12)\\ +16 s_{1} s_{2}^5 (5 s_{2}^2-3)+18 s_{2}^{10}-51 s_{2}^8+92 s_{2}^6-60 s_{2}^4+24 s_{2}^2-4)\pmb{\Bigg]}, 
\end{multline}
where
\begin{equation}
A(s_{1},s_{2})=\pmb{(}9 s_{1}^4+2 s_{1}^2 \left(5 s_{2}^2-3\right)+4 s_{1} s_{2}+s_{2}^4-2 s_{2}^2+2\pmb{)} \pmb{(}s_{1}^4+2 s_{1}^2 \left(5 s_{2}^2-1\right)-4 s_{1} s_{2}+9 s_{2}^4-6 s_{2}^2+2\pmb{)}
\end{equation} 
\end{widetext}

Armed with the Riemannian metric and its corresponding scalar curvature R, we now seek to determine the character of the trajectories in the phase space of Eq.~(\ref{eq:ex01}). Of exceptional interest here is the sign of R.

To discuss the positivity of Eq.~(\ref{eq:R_ex01}), we must demonstrate that the Euler characteristic is positive, i.e. $\chi(\mathcal{M})>0$. Our procedure here is analogous to the one applied to the van der Pol oscillator in Ref.~\cite{Vinicius}. To do that, we shall proceed as follows. First, we determine the nature of the equilibrium point of Eq.~(\ref{eq:ex01}) and its respective index, following the Poincaré-Hopf index theorem. Second, we devote ourselves to the computation of the Euler-characteristic number $\chi(\mathcal{M})$. Finally, the sign of R is interpreted with the aid of the positivity condition expressed by Eq.~(28) of Ref.~\cite{Vinicius}.

To carry out this program, we must observe that Eq.~(\ref{eq:ex01}) has a unique equilibrium point. The Jacobian matrix of Eq.~(\ref{eq:ex01}) at $p_{1}=(0,0)$ is simply 
\begin{equation}\label{eq:JLS} 
\textbf{J}(p_{1})=\begin{bmatrix} 
-1 & -1 \\ 
1 & -1
\end{bmatrix},
\end{equation} 
whence one may conclude that the trace of \textbf{J} is equal to -2 and determinant of $\textbf{J}$ is equal to 2. Owing to the fact that $\mbox{det}\, \textbf{J}>0$ and $\mbox{Tr}\, \textbf{J}^{2}-4\mbox{det}\, \textbf{J}<0$, then $p_{1}=(0,0)$ designates a source equilibrium point. As a final step, we immediately identify, from Theorem 2 of Ref.~\cite{Vinicius}, the index of $p_{1} $ or the Euler characteristic of $\mathcal{M}$ to be   
\begin{equation}\label{eq:xex01} 
\chi(\mathcal{M})= I_{1}(p_{1})=1. 
\end{equation}

Since Eqs.~(28) of Ref.~\cite{Vinicius} and (\ref{eq:xex01}) assure us that $\mbox{R}>0$, the positivity of the scalar curvature is demonstrated. 
 
To conclude, we have studied the sign of scalar curvature of Eq.~(\ref{eq:ex01}) through the geometric methods of Ref.~\cite{Vinicius}. It becomes evident that the scalar curvature R is indeed positive. This result directly contradicts the claims made by Buzzi and Novaes. Their oversight highlights a fundamental error in their analysis, further reinforcing the validity and correctness of our approach.

As a final remark upon the first example\footnote{and, consequently, the third supposed counterexample of Buzzi and Novaes' note, which it is essentially the first example explored above with a transformation of coordinates.} explored above, we must observe that Buzzi and Novaes had omitted another critical information of Eq.~(\ref{eq:ex01}). It should be emphasized that the general form of Eq.~(\ref{eq:ex01}) is actually given by
\begin{equation}\label{eq:ex01a}
\left\{\begin{array}{lcl} &\displaystyle\frac{ds_{1}}{d\tau}&=\Phi_{1}\left(s_{1},s_{2}\right)=-s_{2} + s_{1}(s_{1}^2 + s_{2}^2 - m) ,\\ \\ &\displaystyle\frac{ds_{2}}{d\tau}&= \Phi_{2}\left(s_{1},s_{2}\right)= s_{1} + s_{2}(s_{1}^2 + s_{2}^2 - m) , \end{array}  \right. 
\end{equation}
where $m$ is a control parameter. 

From Eq.~(\ref{eq:ex01a}), we may recognize that it has three varying parameters: $s_{1}$, $s_{2}$, and $m$. Hence, the Riemannian metrical structure of the parameter space of $X=(s_{1},s_{2}, m)$ corresponds, in fact, to a three-dimensional Riemannian geometry. To explore the singularities of the scalar curvature and confirm that R diverges at points symmetric to the origin—indicating the presence of a unique limit cycle—we should note that the curvature R can no longer be calculated using Eq.~(8) from Ref.~\cite{Vinicius}. Instead, the general program to write R is the proper computation of the Christoffel symbols and, then, the fourth-rank curvature tensor, in agreement with the final discussion of Section 03 of Ref.~\cite{Vinicius} and Ref.~\cite{Viniciusprime}. Thus, to  calculate R one must proceed as follows~\cite{Krey}. 
First, one needs to obtain the Christoffel symbols. That is, 
  \begin{equation}\label{Eq:Chris}
\Gamma^{\alpha}_{\xi \, \sigma}=\frac{1}{2}G^{\mu \, \alpha}\left(G_{\mu \,\xi,\,\sigma}+G_{\mu \, \sigma,\,\xi}-G_{\xi \, \sigma,\,\mu}\right), 
 \end{equation}  
in which, the comma notation denotes partial differentiation to the parameters  $s_{1},s_{2}$, and $m$. 
Based on the above, we calculate the curvature tensor by employing
\begin{align}\label{Eq:TensorR}
\text{R}^{\alpha}_{\xi \,\sigma\,l}=\Gamma^{\alpha}_{\xi\,\sigma,\,l}-\Gamma^{\alpha}_{\xi \, l, \, \sigma}+\Gamma^{\mu}_{\xi \,\sigma}\Gamma^{\alpha}_{\mu \,l}-\Gamma^{\mu}_{\xi \, l}\Gamma^{\alpha}_{\mu \,\sigma}. 
\end{align}  
Finally, the Riemannian curvature scalar is 
\begin{align}\label{Eq:CurvRR}
\text{R}=G^{\mu \, \nu}\text{R}^{\vartheta}_{\mu \vartheta \nu}. 
\end{align}

As one can readily see, the derivation of curvature R in the above outlined manner involves nontrivial algebraic computations because the parameter space is a tridimensional geometry. Nevertheless, by computing the scalar curvature R following the outlined program, it can be observed that R diverges at symmetric points with respect to the origin. This behavior indicates that the system exhibits a single limit cycle, in remarkable agreement with the covariant definition of limit cycles of our work. Indeed, it serves as a clear example of a Hopf bifurcation~\cite{Viniciusthird,Viniciusfourth,Viniciusfifth,Viniciussixth}. 

A paper dealing with these details of the GBT and its applications on Hopf bifurcations is now under review.

\subsection{Example 02}
Let us dedicate ourselves to another supposed ``counterexample", carefully manipulated by Buzzi and Novaes. That is to say,

\begin{equation}\label{eq:ex02} 
 \left\{\begin{array}{lcl} &\displaystyle\frac{ds_{1}}{d\tau}&=s_{1} \left(s_{1}^2+s_{2}^2-1\right) \left(s_{1}^2+s_{2}^2-4\right)-s_{2},\\ \\ &\displaystyle\frac{ds_{2}}{d\tau}&=s_{2} \left(s_{1}^2+s_{2}^2-4\right) \left(s_{1}^2+s_{2}^2-1\right)+s_{1},\end{array}  \right. 
\end{equation} 
where $p_{1}=(0,0)$ is an equilibrium point. We note, parenthetically, that Eq.~(\ref{eq:ex02}) is a generalization of Eq.~(\ref{eq:ex01}).

Once more, Buzzi and Novaes state that the scalar curvature of Eq.~(\ref{eq:ex01}) is negative near $p_{1}=(0,0)$. Nevertheless, as we have outlined above and discussed in Ref.~\cite{Vinicius}, the sign of R is determined not by a simple evaluation at equilibrium points but through the Poincaré-Hopf index theorem along with the Gauss-Bonnet theorem.

With this lesson in mind, we now turn to the study of Eq.~(\ref{eq:ex02}) and its scalar curvature with the aid of the geometrical methods of GBT.
\begin{widetext}
It becomes evident from Eq.~(\ref{eq:ex02}) that our second example has, apparently, only two varying parameters: $s_{1}$ and $s_{2}$. Thus, the Riemannian metrical structure of the parameter space $X=(s_{1},s_{2})$ yields a two-dimensional Riemannian geometry. Based on the above, the  positive definite Riemannian metric of the parameter space of Eq.~(\ref{eq:ex02}), with the aid of Eqs.~(5)-(7) of Ref.~\cite{Vinicius}, is given by 
\begin{multline}\label{eq:mex_3B} 
 d\ell^{2}=2\pmb{\Bigg[}\left(2 s_{1} s_{2} \left(s_{1}^2+s_{2}^2-4\right)+2 s_{1} s_{2} \left(s_{1}^2+s_{2}^2-1\right)+1\right)^2+(2 s_{1}^2 \left(s_{1}^2+s_{2}^2-4\right)2 s_{1}^2 \left(s_{1}^2+s_{2}^2-1\right)\\+\left(s_{1}^2+s_{2}^2-4\right) \left(s_{1}^2+s_{2}^2-1\right))^2 \pmb{\Bigg]}(ds_{1})^{2}+2\pmb{\Bigg[}\left(2 s_{1} s_{2} \left(s_{1}^2+s_{2}^2-4\right)+2 s_{1} s_{2} \left(s_{1}^2+s_{2}^2-1\right)-1\right)^2\\+\left(2 s_{2}^2 \left(s_{1}^2+s_{2}^2-4\right)+2 s_{2}^2 \left(s_{1}^2+s_{2}^2-1\right)+\left(s_{1}^2+s_{2}^2-4\right) \left(s_{1}^2+s_{2}^2-1\right)\right)^2 \pmb{\Bigg]}\left(ds_{2}\right)^{2} 
\end{multline}

~~~As we have discussed in Ref.~\cite{Vinicius}, it is implied by Eq.~(\ref{eq:mex_3B}) that this Fisher metric induces a curvature R on the manifold of the parameter space $X$ of the dynamical system of Eq.~(\ref{eq:ex02}). By Eq.~(8) of Ref.~\cite{Vinicius}, we thus have  
\begin{eqnarray}\label{eq:curvRex3B}
\mbox{R}=B(s_{1},s_{2})^{-1}&&\pmb{\Bigg[}\displaystyle
4 ((450 s_{1}^{22}+25 (54 s_{2}^2-359) s_{1}^{20}+(-49122 s_{2}^4-23630 s_{2}^2+76600) s_{1}^{18}-2288 s_{2} s_{1}^{17}\nonumber\\ &&+(-345654 s_{2}^6+689565 s_{2}^4+172776 s_{2}^2-366320) s_{1}^{16}+32 s_{2} (890-213 s_{2}^2) s_{1}^{15}\nonumber\\ &&-4 (3 s_{2}^2 (85401 s_{2}^6-346770 s_{2}^4+340172 s_{2}^2+56720)-270290) s_{1}^{14}-8 s_{2} (548 s_{2}^4-9880 s_{2}^2+18027) s_{1}^{13}\nonumber\\ &&+(2 s_{2}^2 (-848610 s_{2}^8+5142585 s_{2}^6-10265592 s_{2}^4+6634720 s_{2}^2+785372)-2049655) s_{1}^{12}\nonumber\\ &&+8 s_{2} (316 s_{2}^6+8280 s_{2}^4-44140 s_{2}^2+47185) s_{1}^{11}-6 (5 s_{2}^2 (2 (28287 s_{2}^8\nonumber\\ &&-227605 s_{2}^6+677720 s_{2}^4-889040 s_{2}^2+428106) s_{2}^2+73823)-421393) s_{1}^{10}\nonumber\\ &&+8 s_{2} (1960 s_{2}^6-34199 s_{2}^4+91635 s_{2}^2-66461) s_{1}^9+(3 s_{2}^2 (-341604 s_{2}^{12}+3428390 s_{2}^{10}\nonumber\\ &&-13554400 s_{2}^8+26505760 s_{2}^6-25781560 s_{2}^4+9839245 s_{2}^2+688938)-2022680) s_{1}^8\nonumber\\ &&-16 s_{2} (158 s_{2}^{10}+980 s_{2}^8-22225 s_{2}^4+40509 s_{2}^2-23415) s_{1}^7-2 (3 (57609 s_{2}^8\nonumber\\ &&-693540 s_{2}^6+3421864 s_{2}^4-8890400 s_{2}^2+12890780) s_{2}^8\nonumber\\ &&-29682770 s_{2}^6+9009986 s_{2}^4+773120 s_{2}^2-516160) s_{1}^6+8 s_{2} (548 s_{2}^{12}\nonumber\\ &&-8280 s_{2}^{10}+34199 s_{2}^8-44450 s_{2}^6+17070 s_{2}^2-13117) s_{1}^5+(-49122 s_{2}^{18}\nonumber\\ &&+689565 s_{2}^{16}-4082064 s_{2}^{14}+13269440 s_{2}^{12}-25686360 s_{2}^{10}\nonumber\\ &&+29517735 s_{2}^8-18019972 s_{2}^6+3499440 s_{2}^4+998064 s_{2}^2-327635) s_{1}^4\nonumber\\ &&+8 s_{2} ((852 s_{2}^{10}-9880 s_{2}^8+44140 s_{2}^6-91635 s_{2}^4+81018 s_{2}^2-17070) s_{2}^4+1195) s_{1}^3\nonumber\\ &&+2 (675 s_{2}^{20}-11815 s_{2}^{18}+86388 s_{2}^{16}-340320 s_{2}^{14}+785372 s_{2}^{12}-1107345 s_{2}^{10}\nonumber\\ &&+1033407 s_{2}^8-773120 s_{2}^6+499032 s_{2}^4-202335 s_{2}^2+31518) s_{1}^2\nonumber\\ &&+8 s_{2}^3 (286 s_{2}^{14}-3560 s_{2}^{12}+18027 s_{2}^{10}-47185 s_{2}^8+66461 s_{2}^6\nonumber\\ &&-46830 s_{2}^4+13117 s_{2}^2-1195) s_{1}+63036 s_{2}^2+s_{2}^4 (5 (90 s_{2}^{10}-1795 s_{2}^8+15320 s_{2}^6\nonumber\\ &&-73264 s_{2}^4+216232 s_{2}^2-409931) s_{2}^8+2528358 s_{2}^6-2022680 s_{2}^4+1032320 s_{2}^2-327635)-5780)\pmb{\Bigg]},
\end{eqnarray}
in which
\begin{eqnarray}
B(s_{1},s_{2})=&&\pmb{\Bigg[}(25 s_{1}^8+2 (38 s_{2}^2-75) s_{1}^6+(78 s_{2}^4-310 s_{2}^2+265) s_{1}^4+8 s_{2} s_{1}^3+2 (s_{2}^2-3) (2 s_{2}^2-5) (7 s_{2}^2-4) s_{1}^2\nonumber\\&&+4 s_{2} (2 s_{2}^2-5) s_{1}+s_{2}^2 (s_{2}^2-5) (s_{2}^4-5 s_{2}^2+8)+17)^2 (s_{1}^8+2 (14 s_{2}^2-5) s_{1}^6\nonumber\\ &&+(78 s_{2}^4-170 s_{2}^2+33) s_{1}^4-8 s_{2} s_{1}^3+(76 s_{2}^6-310 s_{2}^4+298 s_{2}^2-40) s_{1}^2\nonumber\\ &&+4 s_{2} (5-2 s_{2}^2) s_{1}+5 s_{2}^2 (s_{2}^2-3) (5 (s_{2}^2-3) s_{2}^2+8)+17)^2  \pmb{\Bigg]} 
\end{eqnarray}
\end{widetext}

\newpage

From the knowledge of R, we shall study its sign. Again, the sign of Eq.~(\ref{eq:curvRex3B}) follows from the condition expressed in terms of the Euler-characteristic, Eq.~(28) of Ref.~\cite{Vinicius}. Proceeding in a manner parallel to that of the first example above, we first obtain the nature of the equilibrium points associated with Eq.~(\ref{eq:ex02}) and their respective indices, closely following the index theorem of Poincaré. Second, we devote ourselves to the computation of Euler-characteristic number $\chi(\mathcal{M})$. Finally, we interpret the sign of scalar curvature through the positivity condition. 

It should be noted that Eq.~(\ref{eq:ex02}) possesses a unique equilibrium point $p_{1}=(0,0)$.  From  Eq.~(\ref{eq:ex02}) and $p_{1}=(0,0)$, the Jacobian matrix is directly evaluated as follows 
\begin{equation}\label{eq:JLS} 
\textbf{J}\bm{(}p_{1}\bm{)}=\begin{boldbmatrix} 
4& -1 \\ 
1 & 4
\end{boldbmatrix}.
\end{equation}

From this, it is not difficult to realize that $\mbox{Tr}\, \textbf{J}=8$ and $\mbox{det}\, \textbf{J}=17$. Hence, this immediately teaches us that $p_{1}=(0,0)$ is also a source equilibrium point. From the Poincaré-Hopf index theorem, we identify that the index $I_{1}\bm{(}p_{1} \bm{)}$ or the Euler characteristic of $\mathcal{M}$ to be  1. Because $\chi(\mathcal{M}) >0$, then we have thereby demonstrated the positiveness of the scalar curvature. This last result confirms the consistency of the positivity condition for the sign of the scalar curvature, which contradicts the result found by Buzzi and Novaes.

Before we proceed to the last supposed ``counterexample", we must observe again that Buzzi and Novaes forget to remember or remember to forget that the general form of Eq.~(\ref{eq:ex02}) is actually written as follows
\begin{equation}\label{eq:ex02a} 
 \left\{\begin{array}{lcl} &\displaystyle\frac{ds_{1}}{d\tau}&=s_{1} \left(s_{1}^2+s_{2}^2-m_{1}\right) \left(s_{1}^2+s_{2}^2-m_{2}\right)-s_{2}
,\\ \\ &\displaystyle\frac{ds_{2}}{d\tau}&=s_{2} \left(s_{1}^2+s_{2}^2-m_{1}\right) \left(s_{1}^2+s_{2}^2-m_{2}\right)+s_{1},\end{array}  \right. 
\end{equation} 
where $m_{1}$ and $m_{2}$ are two control parameters.

From Eq.~(\ref{eq:ex02a}), it becomes evident that the actual form of the dynamic system above has, in fact, four varying parameters: $s_{1}$, $s_{2}$, $m_{1}$, and $m_{2}$. Hence, the Riemannian metrical structure of the parameter space of $X$ corresponds to a four-dimensional Riemannian geometry. Consequently, to investigate the singularities of the scalar curvature and verify that R diverges at points symmetric to the origin—indicating the presence of a unique limit cycle—we must remark that the curvature R can no longer be computed using Eq.~(8) from Ref.~\cite{Vinicius}. That is to say, to investigate the singularities of R, one must carry out the program previously mentioned, and construct the curvature tensor\cite{Vinicius,Viniciusprime}.

\subsection{Example 03} 
As the final supposed ``counterexample", we have the following system of differential equations
\begin{eqnarray}\label{eq:ex03} 
\left\{\begin{array}{lcl} &\displaystyle\frac{ds_{1}}{d\tau}&=\Phi_{1}\left(s_{1},s_{2}\right)= -s_{2}+s_{1}^2,\\ \\ \\ 
&\displaystyle\frac{ds_{2}}{d\tau}&= \Phi_{2}\left(s_{1},s_{2}\right)= s_{1}+s_{1}s_{2}.\end{array}  \right. 
\end{eqnarray}

 According to Buzzi and Novaes, the scalar curvature of Eq.~(\ref{eq:ex03}) is positive, and $|\mbox{R}|$ diverges to infinity at $(0,1)$. Therefore, one should expect that Eq.~(\ref{eq:ex03}) exhibits a limit cycle. However, they noted that this last example has only periodic solutions, raising doubts about whether the positivity of $\mbox{R}$ alone is enough to guarantee the existence of limit cycles.

It should be emphasized that Buzzi and Novaes' analysis again neglects two critical aspects of the covariant definition of limit cycles. First, within the general framework of GBT, the sign of R is determined using the Poincaré and Gauss-Bonnet theorems, not just by evaluating equilibrium points. Second is the symmetry of the singularities of $\mbox{R}$. These misinterpretations and the oversight of essential elements of GBT have led to their incorrect conclusions.

With all this before us, we now examine whether it is possible to determine the non-existence of limit cycles of Eq.~(\ref{eq:ex03}) through the geometrical methods of GBT.

We start from the covariant formulation of GBT presented in section 02 of Ref.~\cite{Vinicius}. It is evident from Eq.~(\ref{eq:ex03}) that our third example has only $s_{1}$ and $s_{2}$ as varying parameters. So, the Riemannian metrical structure of the parameter space $X$ corresponds to a two-dimensional Riemannian geometry since Eq.~(\ref{eq:ex03}) has no control parameters.

With the aid of Eqs.~(5)-(7) of Ref.~\cite{Vinicius}  and Eq.~(\ref{eq:ex03}), we see that
\begin{multline}\label{eq:ex2siiimetric} 
 d\ell^{2}=2\pmb{\big[}\left(4 s_{1}^2+(s_{2}+1)^2\right)
\pmb{\big]}^{2}(ds_{1})^{2}\\
 +2\pmb{\big[}\left(s_{1}^2+1\right)\pmb{\big]}\left(ds_{2}\right)^{2}, 
\end{multline}   
which corresponds to an invariant positive definite metric. As a result, the Fisher metric, Eq.~(\ref{eq:ex2siiimetric}), imposes the following curvature R on the two-dimensional parameter space $X$ of Eq.~(\ref{eq:ex03}),

\begin{equation}\label{eq:sIII3Cex3R} 
\mbox{R}=\frac{1}{\left(s_{1}^2+1\right)^2 \left(4 s_{1}^2+(s_{2}+1)^2\right)}.
\end{equation}

 From the knowledge of R, we shall now investigate the trajectories of the phase space of Eq.~(\ref{eq:ex03}) through a general analysis of the sign and the magnitude of scalar curvature R. 

In a manner quite analogous to that presented for the earlier examples, we begin by demonstrating the positivity of R with the aid of Eq.~(28) of Ref.~\cite{Vinicius}. It should be noted that Eq.~(\ref{eq:ex03}) has only one equilibrium point $p_{1}=(0,0)$, which is a center, following the analysis of the Jacobian matrix of Eq.~(\ref{eq:ex03}). From Poincaré-index theorem, we may identify that the index $I\bm{(}p_{1} \bm{)}$ or the Euler characteristic of $\mathcal{M}$ to be
\begin{eqnarray}\label{eq:Xex03} 
\chi(\mathcal{M})=I_{1}\bm{(}p_{1} \bm{)}=1. 
\end{eqnarray}

Inasmuch as $\chi(\mathcal{M})>0$, the defining property of the  Euler-characteristic guarantees the positiveness of Eq.~(\ref{eq:sIII3Cex3R}). That is, Eqs.~(28) of Ref.~\cite{Vinicius} and (\ref{eq:Xex03}) assures us the positivity of the scalar curvature of our last example. From the geometrical interpretation of the sign of curvature R, we thus may conclude that the state space of Eq.~(\ref{eq:ex03}) has periodic trajectories.

With the positiveness of R in mind, we shall now dedicate ourselves to the study of the divergent behavior of $|\mbox{R}|$. By examining the magnitude of scalar curvature, it naturally follows that $|\mbox{R}|$ only diverges to infinity at $(0,1)$. In accordance with the covariant definition of limit cycles~\cite{Vinicius},  on account of the fact that R is positive and it does not diverge at symmetrical points with respect to the origin, then Eq.~(\ref{eq:ex03}) does not possess a limit-cycle, but rather a periodic solution, in agreement with Ref.~\cite{Vinicius}.

In summary, we have investigated the so-called ``counterexamples" presented by Buzzi and Novaes through the geometrical methods of GBT. Our analysis clearly demonstrates that these ``counterexamples" highlight the effectiveness of GBT's and confirm the validity of the covariant definition of limit cycles obtained in our recent work. In addition, after a carefull reading of Buzzi and Novaes' unpublished note, it became evident that  inaccuracies in their conclusions arise from a cursory reading of our work, the omission of critical aspects of the GBT, and reliance on incorrect assumptions or interpretations that deviate from those established in our framework, which led Buzzi and Novaes to erroneous conclusions. This further reinforces the rigor of our approach.

\section{\label{sec:conclusion} CONCLUSIONS}

In the above,  we have investigated the concerns raised by Buzzi and Novaes regarding the Hilbert number and the covariant definition of limit cycles, including a detailed study of the supposed ``counterexamples" they presented.

Initially, we focused on clarifying the apparent contradiction between the Hilbert number presented in our original paper and those proposed earlier. We realize the historical results that suggest a logarithmic growth of $\mbox{H}_{n}$ are fundamentally based on the assumption that the degree of the polynomials of the original sixteenth problem of Hilbert can be modified by a transformation of coordinates with the aid of perturbation theory, which is something that Hilbert never established in the original formulation of the problem.

Moreover, the method introduced by Christopher and Lloyd, and, later, corrected by Li, Chun, and Chan is only devoted to studying the weakened version of Hilbert’s 16th problem posed by V.I.Arnold in 1977, which consists in the study of symmetric planar polynomial Hamiltonian systems and not general polynomial differential equations of degree $n$, as specifically questioned by Hilbert in 1900.

Somewhat surprisingly, when perturbation methods are employed to Hilbert's sixteenth problem without modifying the degree of the polynomials, it has been demonstrated by several authors that $\mbox{H}_{n}$ grows proportionally to $n^{2}$, in agreement with the results obtained by Il'jašenko~\cite{Yas}. In addition, this conclusion is consistent with the findings of Sheng, who demonstrated this result through the use of other geometrical methods. This supports the understanding that the Hilbert number grows as $n^{2}$ and supports the conclusions found in our original work.

Beyond this clarification, we have also investigated the so-called counterexamples of Buzzi and Novaes. Where their calculations went wrong was in assuming that the study of the sign of R is made by a simple evaluation of the scalar curvature at equilibrium points. As well-explained in Ref.~\cite{Vinicius}, the sign of R in the GBT is determined with the aid of the Poincaré and the Gauss-Bonnet theorems. A careful reading of Buzzi and Novaes' note shows that, through relatively trivial errors, they obtained wrong results in virtue of a superficial reading of our work, the omission of critical aspects of the GBT, and reliance on incorrect assumptions or interpretations that diverge from those established in our original work.

With all this before us, we then are led to conclude that the results found in Ref.~\cite{Vinicius}, along with the pieces of evidence presented here, clearly demonstrate the validity of the covariant definition of limit cycles and are in complete agreement with those obtained some time ago in somewhat different circumstances.

The history of science is defined by spectacular debates between opposing ideas, such as Boltzmann's extensive entropy versus Tsallis' non-extensive entropy, the AC versus DC battle between Tesla and Edison, the calculus approaches of Newton and Leibniz, and today, we observe another one: quadratic vs. logarithmic growth for the Hilbert number.

These debates are crucial as they shape the future development of science. Nevertheless,  in today's world, where superficiality is becoming increasingly common, we must acknowledge our responsibility to report results transparently—clearly explaining the conditions under which they were obtained and the hypotheses that underpin them. Otherwise, the omission of these critical details and the inevitable disputes that follow will inevitably lead to the death of science, both nationally and worldwide.

The unpublished note of Buzzi and Novaes\cite{Bossais} is a striking example of this issue. Their omission of the fundamental hypotheses supporting the logarithmic growth result for the Hilbert number and the relevant details of the sacred geometric bifurcation theory was not accidental. On the contrary, it was a deliberate attempt to cast doubt. This behavior hampers progress, disseminates misinformation, and erodes the integrity of scientific and academic discussions. Without honest and well-informed debates, the path ahead becomes unclear, hindering the core of scientific advancement.

\begin{acknowledgments}
V. B. da Silva wishes to express his indebtedness to Yeshua for many enlightening insights. In addition, the authors wish to thank the support and recommendation of professores R. Akhundov, A. Maioli, L. M. Junior, M. H. Francisco, A. L. P. Livorati, L. H. A. Monteiro, J. A. de Oliveira, R. O. M. Torricos, A. Rossi, B. D. Djandue, D. R. da Costa, A. Sahay, N. K. Vitanov, L. A. Barreiro, M. Vukovic as well as C. M. Kuwana, J. D. V. Hermes, F. A. O. Silveira, A. K. Pasquinelli da Fonseca, L. K. A. Miranda, R. M. V. Rocha, and L. H. Pozzo.
\end{acknowledgments}

\newpage

\section*{References}
\nocite{*}
\bibliography{aipsamp}

\end{document}